# HYPOELLIPTICITY: GEOMETRIZATION AND SPECULATION

MICHAEL CHRIST

ABSTRACT. To any finite collection of smooth real vector fields $X_j$ in $\mathbb{R}^n$ we associate a metric in the phase space $T^*\mathbb{R}^n$. The relation between the asymptotic behavior of this metric and hypoellipticity of $\sum X_j^2$, in the smooth, real analytic, and Gevrey categories, is explored.

To Professor P. Lelong, on the occasion of his 85th birthday, in admiration.

## 1. INTRODUCTION

Let $\{X_j\}$ be a collection of real vector fields with $C^\infty$ or $C^\omega$ coefficients, defined in a neighborhood of a point $x_0 \in \mathbb{R}^d$. Consider a second order differential operator $L = \sum_j X_j^2$. Under what conditions is $L$ hypoelliptic, in $C^\infty$ or $C^\omega$? As is well known, certain of these sums of squares operators provide a good model for aspects of the $\bar{\partial}$–Neumann problem in $\mathbb{C}^2$, and a grossly oversimplified but still illuminating window into the higher dimensional case. Indeed, application of the method of boundary reduction to the $\bar{\partial}$–Neumann problem leads to a pseudodifferential equation on the boundary; microlocalizing and composing with an elliptic factor then leads in $\mathbb{C}^2$ to an operator closely related to $X^2 + Y^2$ for certain vector fields. We believe that the study of more general collections of vector fields than those directly relevant to the $\bar{\partial}$–Neumann has shed, and will shed, light on the special cases of interest in complex analysis.

In a companion paper [7] we have established a sufficient condition for the $C^\infty$ case which subsumes a variety of previously known results and examples. This condition involves the existence of auxiliary operators, having certain favorable commutation relations with $L$.

Except in very special cases, there remains a large gap between the hypotheses of the sufficient conditions for hypoellipticity in various function spaces, and counterexamples suggesting what hypotheses are necessary. Many works on various types of hypoellipticity and regularity have relied on commutation methods to achieve results in the positive direction, whereas negative results have been attained through different methods.

How closely is hypoellipticity linked with the commutation method? The aim of this note is to begin to directly investigate this question. We will:

1. Define a metric $\rho_L$, adapted to $L$, on the cotangent bundle, and show how a necessary condition for the existence of auxiliary operators having desirable commutation properties with $L$ may be formulated in terms of $\rho_L$.

---

Research at MSRI is supported in part by NSF grant DMS-9701755. Research also supported by NSF grant DMS-9623007.

2. Speculate on the connection between properties of $\rho_L$ and hypoellipticity of $L$ in $C^\omega$, in Gevrey classes, and to a lesser extent in $C^\infty$.
3. Analyze the behavior of $\rho_L$ in fundamental examples.

A well known construction associates to the collection $\{X_j\}$ of vector fields a Carnot-Caratheodory, or control, metric on the base space $\mathbb{R}^d$, rather than on the cotangent bundle. This metric, and its connection with the fundamental solution of $\sum X_j^2$, has been intensively studied by Fefferman and Phong [12], by Nagel, Stein and Wainger [17], and by Sánchez-Calle [19]. A cotangent bundle metric has been defined by Fefferman and Parmeggiani [11],[18]. Their metric differs from $\rho_L$ in fundamental respects, as explained in §4 below.

## 2. Definitions

Let $\mu : T^*\mathbb{R}^d \mapsto (0, \infty)$ be any positive nonvanishing real symbol belonging to the class $S^1_{1,0}$, and let $\{V_j\}$ be any finite collection of vector fields, defined on a conic open subset of the cotangent bundle $T^*\mathbb{R}^d$, with real $C^\infty$ coefficients. Assume that each $V_j$ takes the form $\sum_k a_{j,k} \partial/\partial x_j + \sum_k b_{j,k} \partial/\partial \xi_j$ where $a_{j,k} \in S^0_{1,0}$ and $b_{j,k} \in S^1_{1,0}$.

Our first definition is made only to facilitate the main definitions below.

**Definition 1.** *A locally $C^1$ function $\psi : T^*\mathbb{R}^d \mapsto \mathbb{R}$ is said to be Lipschitz relative to $(\mu, \{V_j\})$ if at every point $p = (x, \xi)$ of $T^*\mathbb{R}^d$,*

$$(1) \qquad\qquad |V_j\psi(p)| \le \mu(p) \qquad \text{for all } j \le n$$

*and*

$$(2) \qquad\qquad \mu(p)\langle\xi\rangle^{-1}|\nabla_x\psi(p)| + \mu(p)|\nabla_\xi\psi(p)| \le \mu(p).$$

The main condition here is (1), which imposes a stringent upper bound on the variation of $\psi$ in the directions $V_j$, at places where $\mu(x, \xi)$ is smaller than $\langle\xi\rangle$ and some vector field $V_j$ is relatively large. The common factor of $\mu(p)$ on both sides of (2) is not a typographical error and will be explained below.

Next, let finitely many $C^\infty$ real vector fields $X_j$ be given on an open subset of the base space $\mathbb{R}^d$. Recall that to any smooth function $f : T^*\mathbb{R}^d \mapsto \mathbb{R}$ is associated the Hamiltonian vector field

$$H_f = \sum_{n=1}^d \left( \frac{\partial f}{\partial x_n} \frac{\partial}{\partial \xi_n} - \frac{\partial f}{\partial \xi_n} \frac{\partial}{\partial x_n} \right) \ .$$

Define $\sigma_j(x, \xi)$ to be the principal symbol of $iX_j$, which is real valued. Denote by $H_{\sigma_j}$ the Hamilton vector field associated to $\sigma_j$.

In most of the discussion we assume the bracket hypothesis to be satisfied, to some order $m$. Then consider the set $S$ consisting of all multi-indices $I = (I_1, \ldots I_k)$ of degree $|I| = k \in \{1, 2, \ldots m\}$ such that each $I_s \in \{1, 2, \ldots n\}$. For each $I \in S$ define $\sigma_I(x, \xi)$ to be the iterated Poisson bracket $\{\sigma_{I_k}, \sigma_{I'}\}$ where $I' = (I_1, \ldots I_{k-1})$. Define the effective symbol[1] $\tilde{\sigma}$ of $L$ to be

$$(3) \qquad \tilde{\sigma}(x,\xi)^2 = \sum_{I \in S} |\sigma_I(x,\xi)|^{2/|I|}.$$

The bracket hypothesis is equivalent to the existence of $C, \varepsilon \in \mathbb{R}^+$ and $m$ such that

$$\tilde{\sigma}(x,\xi) \geq C\langle\xi\rangle^\varepsilon \qquad \text{for all } (x,\xi).$$

The symplectic structure of the cotangent bundle is reflected in the next definition.

**Definition 2.** *A $C^1$ function $\psi$ defined on $T^*\mathbb{R}^d$ is said to be microlocally Lipschitz relative to $\{X_j\}$ if it is Lipschitz relative to $(\mu, \{V_j\}) = (\tilde{\sigma}, \{H_{\sigma_j}\})$ where $\tilde{\sigma}$ is the effective symbol defined in (3), and $H_{\sigma_j}$ is the Hamiltonian vector field associated to the principal symbol $\sigma_j$ of $iX_j$.*

With this definition, even when $L = \sum X_j^2$ is elliptic, (1) by itself does not force $\psi$ even to be continuous; when $L$ has constant coefficients, for instance, each $V_j = H_{\sigma_j}$ is everywhere in the span of $\{\partial/\partial x_l\}$. The auxiliary condition (2), which does not depend on the vector fields $X_j$ at all because the same factor of $\mu$ occurs on both sides, does ensure that $\psi$ is locally Lipschitz with respect to the standard metric

$$d\rho_0^2 = \langle\xi\rangle^2 dx^2 + d\xi^2.$$

(2) is simply one of the inequalities defining the symbol class $S_{1,0}^1$, rewritten so as to more closely resemble (1). (2) implies (1) when $L$ is elliptic, but not in general.

Our main definition is as follows. Associated to $\{X_j\}$, via $\tilde{\sigma}$ and the collection of Hamiltonian vector fields $\{H_{\sigma_j}\}$, is a metric $\rho$ on $T^*\mathbb{R}^d$.

**Definition 3.** *Let $X_j$ be real $C^\infty$ vector fields in an open subset of $\mathbb{R}^d$. Then*

$$(4) \qquad \rho_{\{X_j\}}(p,q) = \sup_\psi |\psi(p) - \psi(q)|,$$

*where the supremum is taken over all functions $\psi$ that are microlocally Lipschitz relative to $\{X_j\}$. To $L = \sum X_j^2$ is associated the metric*

$$(5) \qquad \rho_L = \rho_{\{X_j\}}.$$

The motivation for this definition stems from a proof of positive hypoellipticity results; see [7].

$\rho_L$ is defined as a metric in the sense of point set topology, but the definition could be modified to give a Riemannian metric whose properties are equivalent for our purpose. (2) implies that for any collection of vector fields $X_j$, there exists $C < \infty$ such that $\rho_{\{X_j\}} \leq C\rho_0$.

One property of $\rho_L$ is that $L$ is elliptic if and only if $\rho_L$ is comparable (on the infinitesimal level) to $\rho_0$.

---

[1]In the absence of the bracket hypothesis, there is no satisfactory general description of the effective symbol. Lemma 3.1 of [7] describes one narrow class of nonsubelliptic operators for which there exists an effective symbol which is significantly larger than the principal symbol.

**Remark.** The construction of $\rho_L$ is conformally invariant in the sense that if $L_2 = r^2 \cdot L_1$ for some constant $r \neq 0$, and if $L_1$ is a subelliptic sum of squares operator, then $\rho_1 \equiv \rho_2$.

If $L_1 = \sum_j X_j^2$ and $L_2 = \sum_j (r(x)X_j)^2$ where $0 \neq r \in C^\infty$, then the two metrics have the same asymptotic behavior, in the sense that their ratio tends to one as the $\xi$ coordinate tends to infinity through any connected conic open set in which $L_1$ has strictly subelliptic behavior in the sense that $C^{-1}\rho_0^\delta \leq \rho_1 \leq C\rho_0^{1-\delta}$ for separated pairs of points.

There exists another construction which roughly speaking is intermediate between $\rho_L$ and the more familiar Carnot-Caratheodory metric, and which will help us to illustrate what phenomenon is measured by $\rho_L$ but not seen by the Carnot-Caratheodory metric. To any finite family $\{X_j\}$ of smooth real vector fields on a connected open set $U \subset \mathbb{R}^d$, this construction associates a one parameter family $\varrho_R$ of metrics on $U$, rather than a single metric on phase space.

Define the auxiliary quantity $\nu : U \times \mathbb{R}^+ \mapsto \mathbb{R}^+$ by

$$(6) \qquad\qquad \nu(x, R) = \min_{|\xi|=R} \tilde{\sigma}(x, \xi).$$

One property of $\nu$ is that for any $R, R' \geq 1$ satisfying $R \sim R'$, one has $\nu(x, R) \sim \nu(x, R')$, uniformly in $x, R, R'$. If the bracket condition holds to order $m$ at $x$, then $\nu(x, R) \geq cR^{1/m}$ as $R \to \infty$.

**Definition 4.** $\varrho_R(p, q)$ *is the supremum of* $|\psi(p) - \psi(q)|$, *taken over all functions* $\psi$ *satisfying*

$$(7) \qquad\qquad |X_j(\psi)(x)| \leq \nu(x, R) \text{ for all } j$$

*and*

$$(8) \qquad\qquad \nu(x, R)R^{-1}|\nabla_x \psi| \leq \nu(x, R).$$

An equivalent formulation is that $\varrho_R$ is the Carnot-Caratheodory metric associated to $\{\nu(x, R)^{-1}X_j\} \cup \{R^{-1}\partial/\partial x_k : 1 \leq k \leq d\}$. If the bracket condition holds to order $m$ at every point of $U$, then $\varrho_R(\cdot, \cdot) \geq cR^{1/m}$ on any compact subset of $U \times U$ minus the diagonal, as $R \to \infty$.

This construction takes a simpler form in the important special case where the set of vector fields $X_j$ has cardinality $d - 1$ and is everywhere linearly independent. Then the characteristic variety of $\sum X_j^2$ is a line bundle. It may be easily checked that a quantity equivalent to $\nu(x, R)$ is obtained by deleting the supremum in (6) and instead taking $\xi$ to be an element of the fiber of that bundle lying over $x$, of magnitude $R$.

## 3. The Question

The following question may be fundamental. Let $s \in [1, \infty)$. We say that two sets $\Gamma, \Gamma' \subset T^*\mathbb{R}^d$ are $c$–separated if any two points $p = (x, \xi) \in \Gamma$, $p' = (x', \xi') \in \Gamma'$ satisfy $|x - x'| + (1 + |\xi| + |\xi'|)^{-1}|\xi - \xi'| > c$. Let $\{X_j\}$ be a finite collection of $C^\omega$ real vector fields.

**Main Question.** *Consider the condition: For every $c > 0$, there exists $c' > 0$ such that for any two $c$–separated open sets $\Gamma, \Gamma' \subset T^*\mathbb{R}^d$,*

$$(9) \qquad \rho_L(p, p') \geq c' \rho_0(p, p')^{1/s} \qquad \text{for all } p \in \Gamma, p' \in \Gamma'.$$

*Is (9) or some closely related condition nearly equivalent to microlocal hypoellipticity of $\sum X_j^2$ in the Gevrey class $G^s$?*

It is essential in (9) that $\Gamma$, $\Gamma'$ be separated; thus no infinitesimal comparison of metrics is made.[2] Restricting attention to pairs $(p, p')$ for which $\rho_0(p, p') \sim 1$, rather than merely $\gtrsim 1$, results in an equivalent question.

In the remainder of this paper, the metric $\rho_L$ will be analyzed for a number of examples, the optimal $C^\infty$, analytic, and Gevrey class hypoellipticity properties of which are already known for each example but the last. In each known case, the result will be seen to be consistent with an affirmative answer to this question, with the word "nearly" omitted.

Nonetheless we have deliberately refrained from posing the question more precisely or labeling it a conjecture. Some variant of (9) may be more directly related to hypoellipticity. Some such variants impose stronger restrictions on $\psi$, involving upper bounds on its partial derivatives of orders greater than one, requiring it to belong to an appropriate symbol class; thus when $s = 1$, $\psi$ should perhaps be an analytic symbol as defined for instance in [21]. Our question merely captures a first order approximation to such a condition. Such variants of (9) are genuinely inequivalent, as shown in Example 9 below. To the extent of this author's knowledge, such a stronger variant is indeed satisfied with $s = 1$ for every operator that has so far been proved to be analytic hypoelliptic.

As will be seen in the subsequent discussion, (9) for $s = 1$ is related to a conjectured criterion for analytic hypoellipticity proposed by Treves [22]. It would be worthwhile to further analyze this connection. Our final Example 9 is an interesting test case which may distinguish between the two points of view.

For $C^\infty$ hypoellipticity the analogous condition would be that for every $A < \infty$ there should exist $C_A < \infty$ such that

$$(10) \qquad \rho_L(p, p') \geq A \log \rho_0(p, p') - C_A$$

for all $p \in \Gamma, p' \in \Gamma'$, where $\rho_L$ would denote the metric associated to an appropriately defined effective symbol of $\sum X_j^2$, if a meaningful effective symbol were to exist.

## 4. Examples

The notation $A \sim B$ will mean that the quantity $A/B$ is bounded above and below by positive finite constants, uniformly in all relevant parameters; usually these parameters will be points in the cotangent bundle or certain subregions thereof, and a large parameter $\lambda$. Throughout the discussion, the zero section is always understood to be excluded from any cotangent bundle.

---

[2]When $s = 1$, the effective symbol may be replaced by $(1 + \sum_j \sigma_j^2)^{1/2}$ in the definition of $\rho_L$ without affecting the validity of (9).

**Example 1.** Let $L = \sum X_k^2$ be an arbitrary sum of squares operator satisfying the bracket hypothesis to order $m$ at every point of an open set $U \subset \mathbb{R}^d$. Then $\tilde{\sigma}(x, \xi) \geq c(1 + |\xi|)^{1/m}$. Therefore in the region where $|\xi| \sim \lambda$, the functions $\lambda^{1/m} x_j$ and $\lambda^{-(m-1)/m} \xi_j$ are microlocally Lipschitz, up to a constant factor, uniformly for large $\lambda$.

From this it follows that

$$(11) \qquad\qquad \rho_L(p, q) \geq c\rho_0(p, q)^{1/m}$$

whenever $\rho_0(p, q) \geq 1$. This is consistent with an affirmative answer to our main question, a theorem of Derridj and Zuily asserting that such operators are Gevrey class hypoelliptic of every order $s \geq m$, and the existence of such operators which are not Gevrey hypoelliptic of any smaller order, such as (13).

On the other hand, for any operator we have $\rho_L \leq C\rho_0$. In particular, when $L$ is elliptic, $\rho_L$ is pointwise comparable to $\rho_0$.

**Remark.** This conclusion for elliptic operators underscores the contrast between $\rho_L$ and a metric $\tilde{\rho}_L$ introduced by Fefferman and Parmeggiani [11],[18]. Firstly, the correspondence between $L$ and $\tilde{\rho}_L$ has a fundamental monotonicity property: For any two operators $L_1$ and $L_2$ whose principal symbols are nonnegative and satisfy $\sigma_1(p) \leq \sigma_2(p)$ for every point $p \in T^*\mathbb{R}^n$, one has $\tilde{\rho}_{L_1}(p, p') \geq \tilde{\rho}_{L_2}(p, p')$ for all $p, p'$. The mapping from $L$ to $\rho_L$ has by design no such monotonicity; this reflects the failure of hypoellipticity of $L_1$ (in various function spaces) to imply hypoellipticity of $L_2$ when $\sigma_1 \leq \sigma_2$, in general. The Carnot-Caratheodory metric associated to $L$ on the base space $\mathbb{R}^d$ likewise shares this monotonic dependence on $L$.

Secondly, up to order of magnitude, $\rho_L$ is *largest* when $L$ is elliptic, whereas $\tilde{\rho}_L$ is then smallest.

**Example 2.** The simplest type of obstruction to (9) in the case $s = 1$ is the possible existence of a nonvanishing vector field $V$ in the span of $\{H_{\sigma_j}\}$, such that some integral curve $\gamma$ of $V$ is contained in the characteristic variety $\Sigma = \{p : \sigma_j(p) = 0 \text{ for all } j\}$. Treves [20] has conjectured that this obstruction should preclude analytic hypoellipticity of $L$.

To see that (9) must fail to hold with $s = 1$ in this situation, let $p = (x, \xi)$ and $p' = (x', \xi')$ be any two distinct points of $\gamma$, and for large $\lambda \in \mathbb{R}^+$ set $\lambda p = (x, \lambda\xi)$ and $\lambda p' = (x', \lambda\xi')$. Consider the integral curve $\gamma_\lambda$ joining $\lambda p$ to $\lambda p'$, parametrized by $s \mapsto \exp(sV)(\lambda p)$; we have $\lambda p' = \exp(c_0 V)(\lambda p)$, where $c_0$ is a positive constant independent of $\lambda$, because the symbols $\sigma_j$ of the underlying vector fields are homogeneous functions of degree one on $T^*\mathbb{R}^d$. Moreover $\gamma_\lambda$ is contained in $\Sigma$ by hypothesis. As is well known, this implies a loss of at least one half of a derivative:

$$\tilde{\sigma}(x, \xi) \leq C(1 + |\xi|)^{1/2} \sim \lambda^{1/2}$$

along $\gamma_\lambda$.

Any microlocally Lipschitz function $\psi$ must satisfy $|V(\psi)| \leq \tilde{\sigma} \leq C\lambda^{1/2}$ along $\gamma_\lambda$. Consequently $|\psi(\lambda p) - \psi(\lambda q)| \leq C\lambda^{1/2}$. Taking the supremum over all such $\psi$ yields

$$\rho_{\{X_j\}}(\lambda p, \lambda p') \leq C\rho_0(\lambda p, \lambda p')^{1/2}$$

as $\lambda \to +\infty$. Consequently (9) fails to be satisfied for $s = 1$.

More precisely, if $\tilde{\sigma}(\lambda q) \leq C\lambda^{1/m}$ for every $q$ belong to $\gamma_\lambda$, then

$$(12) \qquad \rho_L(\lambda p, \lambda p') \leq C\rho_O(\lambda p, \lambda p')^{1/m}.$$

A concrete example is the Baouendi-Goulaouic type operator [1] in $\mathbb{R}^3$:

$$(13) \qquad \partial_x^2 + \partial_y^2 + (x^{m-1}\partial_t)^2$$

Take $p = (0, 0, 0; 0, 0, 1)$ and $p' = (0, 1, 0; 0, 0, 1)$ and $V = \partial_y$.

**Remark.** There are two different ways to formulate the conjecture of Treves mentioned in the discussion of the preceding example. One involves curves $\gamma \subset \Sigma$ whose tangent vectors $\dot{\gamma}$ belong to the span of $\{H_{\sigma_j}\}$ at each point of $\gamma$; the other involves curves for which $\dot{\gamma}$ is orthogonal, with respect to the skew symmetric form defined by the canonical symplectic form $\omega$, to the tangent space to $\Sigma$ at each point of $\gamma$. These two conditions on $\gamma$ are in fact equivalent, under the assumption that the collection of vector fields $\nabla_{x,\xi}\sigma_j$ has constant rank in a neighborhood of $p = \gamma(0)$.

Indeed, $\Sigma$ is by assumption a smooth manifold in a neighborhood of $p$ whose tangent space is the set of all vectors $v \in T_p(T^*\mathbb{R}^d)$ satisfying $v(\sigma_j) = 0$ for all $j$. Denote by $W^\perp$ the orthocomplement, with respect to $\omega$, of any subspace $W$ of $T_p\Sigma$. Since $\omega$ is nondegenerate, $(W^\perp)^\perp = W$ for all subspaces. Define $V$ to be the span of $\{H_{\sigma_j}(p)\}$, a linear subspace of $T_p\Sigma$.

Let $w = \dot{\gamma}(0) \in T_p\Sigma$. Then $w \perp_\omega T_p\Sigma$ if and only if $\omega(v, w) = 0$ for all $v$ such that $v(\sigma_j) = 0$ for all $j$. Since $v(\sigma_j) = \omega(v, H_{\sigma_j})$, the latter condition on $w$ may be restated as $\omega(v, w) = 0$ for all $v \in V^\perp$. Thus $w \perp_\omega T_p\Sigma$ if and only if $w \in (V^\perp)^\perp = V$.

**Remark.** A more general conjecture formulated by Treves in [22] involves the construction of a stratification $\Sigma = \Sigma_0 \supset \Sigma_1 \supset \cdots$, based upon the vanishing not only of the symbols $\sigma_j$, but also of their iterated Poisson brackets, as well as on quantities such as the dimension of the span of their gradients and the rank of the quadratic form defined by the restriction to the tangent space of a stratum of the symplectic form. One considers curves $\gamma \supset \Sigma_k \backslash \Sigma_{k+1}$ whose tangents are orthogonal to $T\Sigma_k$, with respect to $\omega$. According to the preceding remark, this means that $\dot{\gamma}$ belongs to the span of the set of Hamiltonian vector fields associated to the collection of all iterated Poisson brackets of the symbols $\sigma_j$ up to a certain order $n$, depending on $k$. The discussion of the preceding example may be generalized to conclude that the existence of such a curve $\gamma$ implies a certain upper bound on $\rho_L$, depending on $n$.

**Example 3.** Suppose that $X_1, \ldots X_{d-1}$ are linearly independent at each point of $U \subset \mathbb{R}^d$, which throughout this example are assumed to satisfy the bracket hypothesis to order two at every point of $U$. Then $\rho_L(p, q) \geq c\rho_0(p, q)^{1/2}$ for separated points, as we have seen in Example 1. In Example 2 it is shown that this lower bound cannot in general be improved.

If however the characteristic variety $\Sigma = \{(x, \xi) : \sigma_j(x, \xi) = 0 \text{ for all } j\}$ is a symplectic manifold, then each Hamiltonian vector field $H_{\sigma_j}$ is everywhere transverse to $\Sigma$. From this is easily deduced the improvement $\rho_L(p, q) \geq c\rho_0(p, q)^1$ for all separated points. This is consistent with the fundamental theorem of Tartakoff and of Treves, which asserts analytic hypoellipticity in the symplectic situation.

Now consider the one parameter family of metrics $\varrho_R$ on $U$ of Definition 4. $\Sigma$ is a line bundle over $U$, and $\nu(x, R) \sim R^{1/2}$ for all $x, R$. From this and the order two bracket hypothesis it is easily deduced that $\varrho_R(x, y)$ is comparable to $R^{1/2}$ times the associated Carnot-Caratheodory metric at distances greater than $R^{-1/2}$ with respect to the latter metric. In particular, $\varrho_R(x, y) \sim R^{1/2}$ for fixed $x \neq y$. This is so, regardless of whether $\Sigma$ is symplectic. Thus the family $\varrho_R$ fails to detect a property of the vector fields which is fundamental for analytic hypoellipticity.

**Example 4.** Let $W \subset \mathbb{R}^d$ be any connected open set. Let $\Gamma$ be a conic open subset of $T^*W$. Suppose that $L$ is elliptic in $\Gamma$. Then for separated points $p, q$ lying in distinct topological components of $T^*W \backslash \Gamma$, $\rho_L(p, q) \geq c\rho_0(p, q)$. Indeed, any $\psi \in S_{1,0}^1$ whose gradient is supported in $\Gamma$ is microlocally Lipschitz, up to multiplication by a bounded constant factor, and there exist such functions satisfying $|\psi(p) - \psi(q)| \sim \rho_0(p, q)$.

More generally, if $\tilde{\sigma}(x, \xi) \geq C|\xi|^\delta$ for all $(x, \xi) \in \Gamma$, then by the same reasoning $\rho_L(p, q) \geq c\rho_0(p, q)^\delta$ for all $p, q$ satisfying the same hypothesis.

The next example concerns $C^\infty$ hypoellipticity. In the $C^\infty$ case one seeks only to have $\rho(p, p') \gg \log \rho_0(p, p')$ for separated $p, p'$. Under the very weak hypothesis that $L$ has an effective symbol $\tilde{\sigma}(x, \xi)$ which tends to infinity as $|\xi| \to \infty$, this bound is easily obtained whenever $|\xi|$ and $|\xi - \xi'|$ are large, because the function $\psi(x, \eta) = \log(|\eta - \xi|)$ is always locally Lipschitz modulo a constant factor, and an arbitrarily large constant multiple of $\psi$ is Lipschitz relative to $\{X_j\}$ for sufficiently large $|\xi|$, provided that $\tilde{\sigma} \to \infty$.

**Example 5.** In $\mathbb{R}^2$ with coordinates $(x, t)$ consider $L = X^2 + Y^2$ where
$$X = \partial_x, \quad Y = a(x)\partial_t.$$
Suppose that $a \in C^\infty$, and that $a(x) = 0$ if and only if $x = 0$. All such operators are $C^\infty$ hypoelliptic [9], regardless of their degree of degeneracy. Let us see that this is consistent with our point of view.

The characteristic variety is $\Sigma = \{x = \xi = 0\}$ (which is symplectic). Because of the remark made two paragraphs above, the condition $\rho_L(p, q) \gg \log \rho_0(p, q)$ is most likely to fail when $p, q \in \Sigma$ and both have the same $\tau$ coordinate. Thus let $p = (0, 0, 0, \tau)$ and $q = (0, 0, 1, \tau)$ with $\tau$ large and positive. Consider $\psi(x, t, \xi, \tau) = t\tau$. Certainly $|\psi(p) - \psi(q)| = \tau$ is $\gg \log \rho_0(p, q) \sim \log \tau$, and we claim that $\psi$ is Lipschitz. Indeed, the two Hamiltonian vector fields are $\partial_x$ and $a(x)\partial_t - \tau a'(x)\partial_\xi$. The former annihilates $\psi$, while applying the latter to $\psi$ gives $a(x)\tau$, which equals the symbol of $Y$ and hence is acceptably bounded, without taking into account the existence of an effective symbol slightly larger than $|\xi| + |a(x)\tau|$.

**Example 6.** Let $(x, t)$ be coordinates in $\mathbb{R}^2$ and $(x, t; \xi, \tau)$ be coordinates in the cotangent bundle. Consider

$$(14) \qquad\qquad L_1 = \partial_x^2 + (x\partial_t)^2$$

$$(15) \qquad\qquad L_2 = \partial_x^2 + (x\partial_t)^2 + (t\partial_t)^2.$$

The associated effective symbols are up to order of magnitude equal to $\tilde{\sigma}_1 \sim 1 + |\xi| + |x\tau| + |\tau|^{1/2}$, and $\tilde{\sigma}_2 \sim 1 + |\xi| + |x\tau| + |\tau|^{1/2} + |t\tau|$. Thus $\tilde{\sigma}_1 \leq C\tilde{\sigma}_2$.

The Hamiltonian vector fields associated to the vector fields comprising $L_1$ are $\partial_x$ and $x\partial_t - \tau\partial_\xi$. For $L_2$ there is the additional Hamiltonian field $t\partial_t - \tau\partial_\tau$. Thus there arise two competing effects in comparing $\rho_{L_2}$ to $\rho_{L_1}$: $\tilde{\sigma}_2$ is larger, but the presence of an additional Hamiltonian vector field introduces additional constraints on microlocally Lipschitz functions $\psi$.

To analyze the comparison, let $\lambda \in \mathbb{R}^+$ be a large parameter and consider the two points $p = (0, 0; 0, \lambda)$ and $q = (0, 0; 0, 2\lambda)$. The function $\psi(x, t; \xi, \tau) = \tau$ is microlocally Lipschitz with respect to the metric associated to $L_1$, because it belongs to $S^1_{1,0}$ and is annihilated by the two Hamiltonian vector fields. Therefore $\rho_{L_1}(p, q) \geq c\lambda$; it is comparable to the distance between these points in the metric associated to any elliptic operator, even though $L$ fails to be elliptic.

Consider instead $L_2$, and consider the restriction of any microlocally Lipschitz function $\psi$ to the line segment joining $p$ to $q$. The Hamiltonian vector field $t\partial_t - \tau\partial_\tau$ is tangent to this line segment, and $\tilde{\sigma}_2$ has order of magnitude $\lambda^{1/2}$ on the whole segment, so the main condition (1) forces $|\psi(p) - \psi(q)| \leq C\lambda^{1/2}$. Thus $\rho_{L_2}(p, q) \leq C\lambda^{1/2} \ll \rho_{L_1}(p, q)$. Indeed, Métivier [16] has proved that $L_2$ is hypoelliptic in $G^s$ only for $s \geq 2$, whereas $L_1$ is $G^s$ hypoelliptic for all $s \geq 1$ by a special case of a theorem of Grušin [14].

Further analysis shows that $\rho_{L_2}(p, q) \geq c\rho_0(p, q)^{1/2}$ for any points $p, q$. $L_2$ is known to be Gevrey hypoelliptic of order $s$ if and only if $s \geq 2$, so this example is consistent with an affirmative answer to our main question.

So far we have derived upper bounds on $\rho_L$ by examining the variation of microlocally Lipschitz functions $\psi$ along integral curves $\gamma_\lambda \subset \Sigma$ of Hamiltonian vector fields. Our next two examples demonstrate subtler effects due to the influence of one parameter families of integral curves $\gamma_\lambda$ which do not quite lie in the characteristic variety $\Sigma$, but merely tend to $\Sigma$ at a certain rate as $\lambda \to \infty$.

**Example 7.** With coordinates $(x, y, t; \xi, \eta, \tau)$ for $T^*\mathbb{R}^3$, consider

$$(16) \qquad \partial_x^2 + (x^{k-1}\partial_y)^2 + (x^{m-1}\partial_t)^2$$

where $2 \leq k \leq m$. Here

$$\tilde{\sigma} \sim |\xi| + |x^{k-1}\eta| + |x^{m-1}\tau| + |\eta|^{1/k} + |\tau|^{1/m}.$$

Consider the associated Hamiltonian vector field $V = x^{k-1}\partial_y - (k-1)x^{k-2}\eta\partial_\xi$. Let $p = (\delta, 0, 0; 0, 0, \lambda)$ where $\delta = \lambda^{-1/m}$, and $q = \exp(TV)(p)$ where $T = \delta^{1-k}$. The integral curve $\gamma$ joining $p$ to $q$ takes the simple form $\exp(sV)(p) = (\delta, s\delta^{k-1}, 0; 0, 0, \lambda)$. In particular, $q = (\delta, 1, 0; 0, 0, \lambda)$. Along $\gamma$, $\tilde{\sigma} \sim \lambda^{1/m}$ because $\delta$ was chosen so that $\delta^{k-1}\lambda \sim \lambda^{1/m}$. Thus any microlocally Lipschitz function $\psi$ must satisfy $|V\psi| \leq C\lambda^{1/m}$, whence $|\psi(p) - \psi(q)| \leq CT\lambda^{1/m} = C\lambda^{(k-1)/m}\lambda^{1/m} = \lambda^{k/m}$.

On the other hand, the function $\psi(x, y, t; \xi, \eta, \tau) = \lambda^{k/m}y$ is microlocally Lipschitz in the region where $|(\xi, \eta, \tau)| \sim \lambda$, and $|\psi(p) - \psi(q)|$ is comparable to $\lambda^{k/m}$. Thus

$$\rho_L(p, q) \sim \lambda^{k/m} \sim \rho_0(p, q)^{k/m}.$$

This is consistent with both an affirmative answer to our main question, and the known fact that $L$ is hypoelliptic in the Gevrey class of order $s$ if and only if $s \geq m/k$, assuming that $m \geq k$ [4].

**Example 8.** Fix any integers $m, r > 1$. In $\mathbb{R}^2$ consider

$$\partial_x^2 + (x^{m-1}\partial_t)^2 + (t^r\partial_t)^2; \tag{17}$$

the case $r = 1$ has been discussed above. One of the associated Hamiltonian vector fields is $V = t^r\partial_t - rt^{r-1}\tau\partial_\tau$. The effective symbol has order of magnitude $|\xi| + |x\tau| + |t^r\tau| + |\tau|^{1/m}$.

Given a large quantity $\lambda$, define $\delta = \lambda^{-(m-1)/mr}$, so that $\delta^r\lambda = \lambda^{1/m}$. Let $p = (0, \delta; 0, \lambda)$ and $q = \exp(-TV)(p)$ where $T = \delta^{1-r}$. Denote by $\gamma$ the segment of integral curve $\exp(-sV)(p)$, $0 \leq s \leq T$. For every point of $\gamma$, uniformly in $\lambda$, $\exp(-sV)(p) = (0, t(s); 0, \tau(s))$ where $\delta \geq t(s) \geq c\delta$ and $\lambda \leq \tau(s) \leq C\lambda$. Moreover $q = (0, t(T); 0, \tau(T))$ where $\tau(T) - \lambda$ is positive and has order of magnitude $\lambda$. By the choice of $\delta$ and the fact that $|t(s)| \leq \delta$ for every point of $\gamma$, $\tilde{\sigma} \leq C\lambda^{1/m}$ at every point of $\gamma$.

Every microlocally Lipschitz function $\psi$ must therefore satisfy $|V\psi| \leq C\lambda^{1/m}$ at every point of $\gamma$, and consequently $|\psi(p) - \psi(q)| \leq CT\lambda^{1/m}$. Since $T = \delta^{1-r}$, we conclude after a small computation that

$$\rho(p, q) \leq C\lambda^{(mr-m+1)/mr}. \tag{18}$$

This upper bound on $\rho_L(p, q)$ is also a lower bound, up to a constant factor. To see this consider the function $\psi(x, t; \xi, \tau) = \lambda^a \cdot \tau$, where $a + 1 = (mr - m + 1)/mr < 1$. Then in the region where $|(\xi, \tau)| \sim \lambda$, $\psi$ is a symbol in $S_{1,0}^1$ uniformly in $\lambda$, and $|V(\psi)| \leq C|t|^{r-1}\lambda^{a+1}$. To show that $\psi$ is in this region microlocally Lipschitz relative to the vector fields in question, up to a bounded factor, we must verify that $|V\psi| \leq C\tilde{\sigma}$. Indeed, this holds because $\tilde{\sigma}(x, t; \xi, \tau) \geq c\lambda^{1/m} + c|t|^r\lambda$ there, and $|t|^{r-1}\lambda^{a+1}$ equals the logarithmically convex average $(|t|^r\lambda)^\theta \cdot (\lambda^{1/m})^{1-\theta} \leq C\tilde{\sigma}$ where $\theta$ is defined to be $(r-1)/r$. Thus $\psi$ is microlocally Lipschitz, and

$$|\psi(p) - \psi(q)| \sim \lambda^{(mr-m+1)/mr}, \tag{19}$$

whence $\rho_L(p, q) \geq c\lambda^{(mr-m+1)/mr}$. It has indeed been proved in [6], as well as in [2] and [15], that $L$ is $G^s$ hypoelliptic for all $s$ satisfying $s^{-1} \leq (mr - m + 1)/mr$. In unpublished work this author has outlined a proof that the threshold $mr/(mr-m+1)$ is optimal, based on the method introduced in [3] and [5].

**Remark.** Both a conjectured characterization of analytic hypoellipticity by Treves [22], and a conjectured sufficient condition for Gevrey hypoellipticity by Bove and Tartakoff [2] are formulated in terms of properties of a certain stratification of the characteristic variety $\Sigma$ of $L$. This stratification does not distinguish between different parameters $r$ in Example 8; $r$ does influence the optimal degree of Gevrey class hypoellipticity.

## 5. A Subtler Example

We next discuss an operator for which analytic hypoellipticity apparently remains an open question. It is the principal part of the Kohn Laplacian for a three dimensional pseudoconvex CR manifold, and is of particular interest both because it is one of the simplest operators whose analytic hypoellipticity is undecided, and because it illustrates the distinction between two variants of our main question; functions $\psi$ which are merely Lipschitz relative to $\rho_L$ can vary much more rapidly than those which also belong to the symbol class $S^1_{1,0}$.

Consider the domain $\Im(z_2) > b(z_1)$ in $\mathbb{C}^2$, where $b$ is assumed to be subharmonic so that this domain will be pseudoconvex. Its boundary may be identified with $\mathbb{R}^3$, with coordinates $(x, y, t)$, in such a way that a Cauchy-Riemann operator is $\bar\partial_b = X + iY$ where $X = \partial_x - b_y \partial_t$ and $Y = \partial_y + b_x \partial_t$. Here $b_x = \partial_x b$, $b_y = \partial_y b$, and $b(x, y) \equiv b(x + iy)$. The Kohn Laplacian is then $(X + iY)(X - iY)$, which equals $X^2 + Y^2$ modulo a lower order term.[3] Analytic or Gevrey hypoellipticity of such an operator depends only on $\Delta b$, where $\Delta$ denotes the Laplacian in $\mathbb{R}^2$, rather than on $b$ itself, for any operator may be transformed into any other with the same invariant $\Delta b$ via a change of variables $(x, y, t) \mapsto (x, y, t - \phi(x, y))$.

**Example 9.** In $\mathbb{R}^3$ consider $L = X^2 + Y^2$ with

$$X = \partial_x - b_y \partial_t, \quad Y = \partial_y + b_x \partial_t$$

where $b$ is a polynomial depending only on $x, y$, satisfying

$$\Delta b(x, y) = x^k + y^k + x^2 y^2,$$

and $k \geq 6$ is an integer. Define $\lambda(x, y) = \Delta b = x^k + y^k + x^2 y^2$.

For general $b(x, y)$ it has been shown by Grigis and Sjöstrand [13] that if $\Delta b$ vanishes only at $x = y = 0$ and is a polyhomogeneous function of $(x, y)$, that is, $\Delta b(r^k x, r^l y) \equiv r^n \Delta b(x, y)$ for some positive integers $k, l, n$, then $L$ is analytic hypoelliptic. In our example $\lambda = \Delta b$ is polyhomogeneous for $k = 4$, but not for larger $k$.

Since $[X, Y] = \lambda(x, y)\partial_t$ and $\lambda$ vanishes only on the line $\{(x, y, t) : x = y = 0\}$ in $\mathbb{R}^3$, the bracket hypothesis is satisfied to order 2 everywhere except on that line. Where $x = y = 0$, it is satisfied to order precisely $m = 6$, because of the presence of the term $x^2 y^2$ in $\lambda$. The characteristic variety $\Sigma \subset T^*\mathbb{R}^3$ of $L$ is a manifold of codimension two. It is symplectic, at every point at which the principal symbol of $[X, Y]$ is nonzero. Thus $\Sigma$ is symplectic everywhere except on the submanifold $\Sigma' \subset \Sigma$ given by $\Sigma' = \{(x, y, t; \xi, \eta, \tau) : x = y = \xi = \eta = 0\}$. Hence by virtue of the theorem of Tartakoff and of Treves, $L$ is microlocally analytic hypoelliptic everywhere except possibly on $\Sigma'$.

Any particular iterated commutator of $X, Y$ has a principal symbol which either vanishes identically on $\Sigma'$, or vanishes nowhere on $\Sigma'$. Therefore the stratification defined by Treves [22] takes the simple form $\Sigma \supset \Sigma' \supset \emptyset$. Because $\Sigma'$ is itself a symplectic manifold, his conjecture then predicts that $L$ should be analytic hypoelliptic.

---

[3] This lower order term is not insignificant, but the questions of interest regarding $(X+iY)(X-iY)$ concern its microlocal hypoellipticity properties in a subset of the cotangent bundle, where its properties coincide with those of $X^2 + Y^2$ in all known cases.

For simplicity we set $k = 6$ for the remainder of the discussion. All results generalize to $k \geq 6$, but the details of the analysis change slightly.

We will establish two contrasting conclusions concerning this example. Set $p_T = (0, 0, 0; 0, 0, T)$ and $q_T = (0, 0, 1; 0, 0, T)$, where $T \in \mathbb{R}^+$ is a large parameter.

**Proposition 5.1.** *There exists $c > 0$ such that for each large $T$ there exists a function $\varphi_T$ which is microlocally Lipschitz relative to $\{X, Y\}$, and satisfies $|\varphi(q_T) - \varphi(p_T)| \geq cT$.*

**Proposition 5.2.** *It is not the case that there exist $c > 0$ and a one parameter family of functions $\varphi_T$ such that for each large $T \in \mathbb{R}^+$, $\varphi_T$ is microlocally Lipschitz relative to $\{X, Y\}$, $|\varphi_T(q_T) - \varphi_T(p_T)| \geq cT$, and moreover*

$$|\partial_{x,y,t}^\alpha \partial_{\xi,\eta,\tau}^\beta \varphi_T| \leq C_{\alpha,\beta}(1 + |(\xi, \eta, \tau)|)^{1-|\beta|} \qquad \text{for all } |\alpha| + |\beta| \leq 5.$$

Thus whether our philosophy predicts analytic hypoellipticity or not, depends on precisely which variant of $\rho_L$ we consider. This author ventures no conjecture as to whether this particular operator is analytic hypoelliptic.

We first prove Proposition 5.1. We will construct a single function $\varphi$ which has the required properties simultaneously for all $T$. Consider a function

$$\varphi(x, y, t; \xi, \eta, \tau) = (t + h)\tau + f \cdot (\xi - \tau b_y) + g \cdot (\eta + \tau b_x)$$

where $f, g, h$ are functions of $x, y$ alone. Our aim is construct $f, g, h$ so that $\varphi$ is microlocally Lipschitz relative to $\{X, Y\}$. Since $\varphi(q) - \varphi(p) = T - 0 = T$, we can then conclude that $\rho_L(p, q) \geq cT$.

Because the characteristic variety of $L$ is $\{x = y = \xi = \eta = 0\}$ and because the function $(x, y, t; \xi, \eta, \tau) \mapsto \tau$ is annihilated by both $H_X, H_Y$, it is easy to conclude using the same function $\varphi$ that more generally, for any separated points $p, q$, one has $\rho_L(p, q) \geq c\rho_0(p, q)$.

The Hamiltonian vector fields associated to $X, Y$ are

$$H_X = \partial_x - b_y \partial_t + \tau b_{yx} \partial_\xi + \tau b_{yy} \partial_\eta,$$
$$H_Y = \partial_y + b_x \partial_t - \tau b_{xx} \partial_\xi - \tau b_{xy} \partial_\eta.$$

Thus

$$H_X \varphi = -\tau b_y - \tau f b_{xy} + \tau g b_{xx} + \tau f b_{xy} + \tau g b_{yy} + \tau h_x \qquad \text{plus error terms}$$
$$= \tau[-b_y + g\Delta b + h_x] \qquad \text{plus error terms}$$
$$H_Y \varphi = -\tau f b_{yy} + \tau g b_{xy} + \tau b_x - \tau b_{xx} f - \tau b_{xy} g + \tau h_y \qquad \text{plus error terms}$$
$$= \tau[b_x - f\Delta b + h_y] \qquad \text{plus error terms.}$$

Each of the aforementioned error terms is a product of one of $f_x, f_y, g_x, g_y$ with either $(\xi - \tau b_y)$ or $(\eta + \tau b_x)$. The latter two quantities are the principal symbols of $X, Y$, hence are dominated by $\tilde{\sigma}$. If $f, g$ are constructed so as to be Lipschitz continuous functions of $(x, y)$, then each error term will be majorized by $\tilde{\sigma}$. Since our goal is to have $|H_X\varphi| + |H_Y\varphi| \leq C\tilde\sigma$, it therefore suffices to construct Lipschitz functions $f, g$ and a function $h$ such that

$$(20) \qquad b_y - g\Delta b = h_x \qquad \text{and} \qquad -b_x + f\Delta b = h_y.$$

Recalling that $\Delta b = \lambda$, a necessary condition for (20) is that $b_{yy} - (g\lambda)_y = -b_{xx} + (f\lambda)_x$, which may be rewritten as

$$(21) \qquad \lambda = (f\lambda)_x + (g\lambda)_y.$$

Conversely, if $f, g$ are bounded solutions of (21) then there exists a solution $h$ of (20) in the sense of distributions. But the left hand sides of both equations in (20) are bounded since $b$ is a polynomial, so $h$ is Lipschitz after possible redefinition on a set of measure zero; in fact $h \in C^{1,1}$ when $f, g$ are Lipschitz. Thus in order to show that $\varphi$ is microlocally Lipschitz relative to $\{X, Y\}$, it suffices to prove the existence of Lipschitz continuous (in the ordinary sense) solutions $f, g$ of (21).

**Lemma 5.3.** *There exist Lipschitz continuous functions $f, g$, defined in a neighborhood of the origin in $\mathbb{R}^2$, satisfying $\lambda \equiv (f\lambda)_x + (g\lambda)_y$.*

The conclusion is to be interpreted as equality almost everywhere. Note that if $\lambda$ were homogeneous or polyhomogeneous, then there would be a solution of the form $f = k_1 x$, $g = k_2 y$ for certain constants $k_j$, by Euler's identity. It is a relatively simple matter to construct Lipschitz solutions of $\lambda = f \cdot \lambda_x + g \cdot \lambda_y$ in the present case, but that is not what is required.

*Proof.* We write $A \lesssim B$ to mean that the ratio $A/B$ of nonnegative functions is bounded above by a finite constant, and $A \sim B$ to mean that $A \lesssim B$ and $B \lesssim A$. Different constructions will be used in different regions of the plane. Consider first the region $\Gamma_1$ where $|x| \leq 2|y|$. We set $g \equiv 0$ and solve $(f\lambda)_x = \lambda$ by defining

$$f(x, y) = \lambda(x, y)^{-1} \int_0^x \lambda(s, y)\, ds.$$

In this region $\lambda \sim y^6 + x^2 y^2$, so $|\int_0^x \lambda(x, y)\, ds| \lesssim |x|\lambda(x, y)$. Hence $|f(x, y)| \lesssim |x|$.

In analyzing $\nabla f$ it will be useful to note that in $\Gamma_1$,

$$|\lambda_x| \lesssim \lambda/|x| \qquad \text{and} \qquad |\lambda_y| \lesssim \lambda/|y|.$$

To analyze $\nabla f$ consider first

$$f_x = 1 - \lambda^{-2}\lambda_x \cdot \int_0^x \lambda(s, y)\, ds.$$

Since $|\lambda_x| \lesssim \lambda/|x|$ and the integral is $\lesssim |x|\lambda$, the second term is uniformly bounded. Similarly

$$f_y = -\lambda^{-2}\lambda_y \int_0^x \lambda + \lambda^{-1}\int_0^x \lambda_y(s, y)\, ds.$$

In absolute value the first term is $\lesssim \lambda^{-2} \cdot (\lambda/|y|) \cdot (|x|\lambda)$, which is uniformly bounded in $\Gamma_1$. The absolute value of the second term is $\lesssim \lambda^{-1}\int_0^x \lambda/|y|\, ds \lesssim |x|/|y|$, so likewise is bounded.

Fix an auxiliary function $\phi$ which is homogeneous of degree zero in $\mathbb{R}^2$, is $C^\infty$ except at the origin, is identically equal to one where $|x| \leq |y|/2$, and is supported where $|x| \leq |y|$. Set $\tilde{f} = f \cdot \phi$. This formula makes sense only where $|x| \leq |y|$, but we extend the definition by setting $\tilde{f} \equiv 0$ where $|x| > |y|$. Then $\tilde{f}$ is Lipschitz, because $|\nabla \phi(z)| \lesssim |z|^{-1}$ and $|f(z)| \lesssim |z|$. We claim that where $|y|/2 \leq |x| \leq |y|$ and $|(x, y)|$ is sufficiently small,

$$|\nabla^\alpha \tilde{f}(z)| \leq C_\alpha |z|^{1-|\alpha|}$$

for every multi-index $\alpha$. Because $\phi$ is homogeneous and smooth except at the origin, it suffices to verify this for $f$. One has explicitly

$$f(x, y) = \frac{xy^6 + x^3y^2/3 + x^7/7}{y^6 + x^2y^2 + x^6}.$$

In the conic region in question, $\lambda(z) \sim |z|^4$, the terms $x^6, y^6$ in the denominator being comparatively negligible. Hence $|\nabla^\alpha f(z)| \lesssim |z|^{1-|\alpha|}$ where $|x| \sim |y|$. The claim follows.

In the (overlapping) region where $|y| \leq 2|x|$, we solve instead $(g\lambda)_y = \lambda$, taking $f \equiv 0$. Because $\lambda$ is a symmetric function of $(x, y)$, conclusions parallel to those above may be obtained. Define $\tilde{g}(x, y) = g(x, y)\phi(y, x)$, and again define $\tilde{g} \equiv 0$ where $|y| > |x|$. We have then $\lambda \equiv (\tilde{f}\lambda)_x + (\tilde{g}\lambda)_y$ where $|x| \leq |y|/2$ and also where $|y| \leq |x|/2$.

Define

$$\tilde{\lambda} = \lambda - (\tilde{f}\lambda)_x - (\tilde{g}\lambda)_y.$$

$\tilde{\lambda}$ is supported in the conic region $\Gamma_3$ where $|y|/2 \leq |x| \leq 2|y|$. It satisfies $|\nabla^\alpha \tilde{\lambda}(z)| \lesssim |z|^{4-|\alpha|}$ for all $z \in \Gamma_3$ and all $\alpha$.

In order to complete the proof of the lemma, it suffices to construct Lipschitz functions $F, G$ satisfying

$$\tilde{\lambda} = (F\lambda)_x + (G\lambda)_y.$$

To accomplish this we set $F = xh$ and $G = yh$ and solve $\tilde{\lambda} = (xh\lambda)_x + (yh\lambda)_y$ for the single unknown $h$. In polar coordinates $(r, \theta)$ the equation becomes

$$\tilde{\lambda} = \lambda r h_r + h \cdot (2\lambda + x\lambda_x + y\lambda_y).$$

Define

$$\beta(x, y) = (2 + x\lambda_x/\lambda + y\lambda_y/\lambda)$$

and note that $2 \leq \beta$ and moreover that in $\Gamma_3$, $|\nabla^\alpha \beta(z)| \lesssim |z|^{-|\alpha|}$. The equation to be solved is

$$(22) \qquad\qquad h_r + h\beta/r = r^{-1}\tilde{\lambda}/\lambda.$$

It suffices to produce a bounded solution whose gradient is $\lesssim r^{-1}$.

Note for future use that

$$|\beta_\theta(r, \theta)| \lesssim r$$

in $\Gamma_3$, not merely $\lesssim 1$. The fact that $\lambda$ is homogeneous of degree 4 modulo higher order terms and nonvanishing in $\Gamma_3$ leads by Euler's identity to the conclusions $\beta(r, \theta) = 6 + O(r)$ and $\nabla\beta = O(1)$ in that cone. This strengthened bound is essential to our construction.

Define next

$$b(r, \theta) = \int_1^r \beta(\rho, \theta), d\rho/\rho.$$

Then $b(r, \theta) < 0$ for $r < 1$, and $|b(r, \theta)| \sim \log r^{-1}$. Because $\beta/\rho \geq 2/\rho$, $e^b(\rho, \theta) \lesssim \rho^2$ as $\rho \to 0$.

Define a solution $h$ of (22) by

$$h(r, \theta) = e^{-b} \int_0^r e^b(\rho, \theta) \frac{\tilde{\lambda}}{\lambda}(\rho, \theta) \, \frac{d\rho}{\rho};$$

the integral converges since $\tilde{\lambda}/\lambda$ is bounded and $e^b \lesssim \rho^2$. Formally $h$ satisfies the required equation. Note that $h$ is supported in $\Gamma_3$, since the integrand vanishes identically outside it.

To complete the proof of the lemma it suffices to show that $|h| \lesssim 1$, and that $|\nabla h(r, \theta)| \lesssim r^{-1}$. Because we have seen above that $\tilde{\lambda}/\lambda$ is bounded, an upper bound for $|h|$ is $e^{-b} \int_0^r e^b \, d\rho/\rho$. Since $b_r$ is between $2/r$ and $C/r$, the factor $e^b$ is monotone increasing, and $e^b(r/2, \theta) \leq \frac{1}{4} e^b(r, \theta)$. Consequently this last integral has the same order of magnitude as $e^b$, and hence $h$ is uniformly bounded.

From the boundedness of $h$ and the differential equation (22) it follows immediately that the partial derivative $h_r \lesssim r^{-1}$ in absolute value. It remains to show that $|h_\theta| \lesssim 1$. Differentiation in the definition of $h$ yields several terms; the simplest is $-b_\theta$ times $h$. We have

$$|b_\theta| \leq \int_r^1 |\beta_\theta| \, d\rho/\rho,$$

and $|\beta_\theta(\rho, \theta)| \leq \rho|\nabla\beta| \lesssim \rho$ in $\Gamma_3$, as observed above. Thus $b_\theta$ is uniformly bounded. Hence this simplest term is $\lesssim h$, hence uniformly bounded.

A second term arises when on differentiating $h$ with respect to $\theta$, the derivative falls upon the factor of $e^b$ inside the integral. An additional factor of $b_\theta$ results, causing no harm since it is uniformly bounded; the analysis indicated above for $h$ itself applies also to this term. The third and last term arises when the derivative falls on the factor of $\tilde{\lambda}/\lambda$ inside the integral. Now

$$|\partial_\theta(\tilde{\lambda}/\lambda)(r, \theta)| \leq Cr|\nabla(\tilde{\lambda}/\lambda)| \lesssim 1$$

in $\Gamma_3$, since $|\nabla\tilde{\lambda}| \lesssim |z|^{4-|\alpha|}$ in $\Gamma_3$, the same holds for $\lambda$, and $|z|^4 \lesssim \lambda(z)$ there. Thus the same analysis applies once more. Hence $|\nabla h| \lesssim r^{-1}$, and the proofs of both the lemma and Proposition 5.1 are complete. $\qquad\square$

*Proof of Proposition 5.2.* Suppose to the contrary that functions $\varphi_T$ possessing the indicated properties were to exist. Note that $H_X, H_Y$ are homogeneous of degree zero with respect to the dilations $(x, y, t; \xi, \eta, \tau) \mapsto (x, y, t; r\xi, r\eta, r\tau)$ for $r \in \mathbb{R}^+$. By applying Arzela-Ascoli and a rescaling argument based on these dilations and homogeneity to the collection $\{\varphi_T\}$, we could then extract a single function $\varphi$ which in any fixed relatively compact subset of $T^*\mathbb{R}^3$ minus the zero section belonged to $C^5$ and satisfied $\varphi(p) \neq \varphi(q)$, where $p = p_1, q = q_1$, and

$$(23) \qquad |H_X\varphi| + |H_Y\varphi| \leq |\xi - \tau b_y| + |\eta + \tau b_x|.$$

The vector fields $H_X, H_Y$ are tangent to each level set of $\tau$, so the restriction of $\varphi$ to $\{\tau = 1\}$ must satisfy the preceding inequality. Henceforth we restrict attention to that level set, freezing $\tau = 1$. Since $\varphi(p) \neq \varphi(q)$, there must exist $t$ such that $\partial_t\varphi(0, 0, t; 0, 0, 1) \neq 0$. Both $H_X, H_Y$ are translation invariant with respect to $t$, so we may translate the coordinates so that $\partial_t\varphi(0, 0, 0; 0, 0, 1) \neq 0$.

Next, expanding in Taylor series with respect to $t, \xi, \eta$ about the origin and rewriting the result slightly gives

$$\varphi(x, y, t; \xi, \eta, 1) = h + \gamma t + f \cdot [\xi - b_y] + g \cdot [\eta + b_x] + E$$

where $E(x, y, t; \xi, \eta) = O(t, \xi, \eta)^2$ and the coefficients $\gamma, f, g, h$ are functions of $x, y$ alone. Moreover $\gamma(0, 0) \neq 0$ because $\partial_t\varphi \neq 0$; by dividing through by a constant we may assume that $\gamma(0, 0) = 1$.

Applying $H_X, H_Y$, evaluating at $t = \xi = \eta = 0$, and invoking (23) yields now the equations (20) that we arrived at in the proof of Proposition 5.1. Consequently $\lambda \equiv (f\lambda)_x + (g\lambda)_y$.

Consider the Taylor expansions of both sides of this last equation about the origin, and compare terms of equal degrees. By examining terms homogeneous of degree three, one finds immediately that $f(0) = g(0) = 0$. Equality for degree 4 forces $f = c_1 x + O^2(x, y)$ and $g = c_2 y + O^2(x, y)$, where $3c_1 + 3c_2 = 1$, in order to reproduce the term $x^2 y^2$ on the left without introducing other monomials of degree 4. The contribution to $(f\lambda)_x + (g\lambda_y)$ of the sum of all terms homogeneous of degree 2 in the Taylor expansions of $f, g$ is a homogeneous polynomial of degree 5 plus a remainder which is $O^7(x, y)$. The part of degree 5 must vanish identically, since no such terms are present in the expansion of $\lambda$. Consider finally the coefficient of $x^6$ on the right. It equals $7c_1 + c_2$, for the monomials $x^2 y^2$ and $y^6$ in $\lambda$ cannot possibly lead to any $x^6$ term on the right. Similarly the coefficient of $y^6$ equals $c_1 + 7c_2$. Therefore $7c_1 + c_2 = 1 = c_1 + 7c_2$. The unique solution is $c_1 = c_2 = 1/8$, but this is incompatible with $3c_1 + 3c_2 = 1$. This completes the proof of nonexistence of $\varphi$.  □

**Summary.** Although a conjecture linking the growth of $\rho_L$ to hypoellipticity, in various function spaces, is consistent with a wide variety of examples, Example 9 suggests that this conjecture may be overly simplistic. Another reason for wariness is that such a conjecture takes no account of the behavior of the symbols $\sigma_{X_j}$ in the complexified phase space. To decide whether Example 9 is analytic hypoelliptic might be quite illuminating.

**Acknowledgement.** The author is indebted to Chun Li for a suggestion used in the proof of Lemma 5.3, and to Charles Fefferman for kindly explaining [11] and [18].

MICHAEL CHRIST, DEPARTMENT OF MATHEMATICS, UNIVERSITY OF CALIFORNIA, BERKELEY 94720, BERKELEY, CA USA

*E-mail address*: mchrist@math.berkeley.edu